\documentclass[a4paper]{article}
\pdfoutput=1

\usepackage{IEK10} % defined in IEK10.sty
\usepackage{natbib}

\makeatother
\usepackage{xcolor}
\usepackage{mathtools, cuted}
\usepackage{tabularx}
\usepackage{eurosym}
\usepackage[utf8]{inputenc}

\usepackage{pstricks, pst-plot}	% pstricks
\usepackage{graphicx} % enables import of various formats
\usepackage{wrapfig} % enables floating wrapped figures
\usepackage[figuresright]{rotating}
\usepackage{lipsum}
\usepackage{placeins}
\usepackage{etoolbox}
\usepackage{amsmath}
\numberwithin{equation}{section} % changes equation numbering to (section.equtionnumber)

\usepackage{setspace}
\usepackage{enumitem}
\usepackage[english]{babel}
\usepackage{blindtext}
\usepackage{epsfig} %hilft evtl mit eps dateien
\usepackage[utf8]{inputenc}			% Verwendung von Utf8-Codierungen -> Bei Darstellungsproblemen 'utf8x' verwenden
\usepackage{ifpdf} 							% Änderungen bei PDF-Erzeugung aktivieren (Bessere Bilddarstellung, etc.)
\usepackage{multicol}						% Ermöglichung von mehrspaltigen Dokumenteilen
\usepackage{titletoc}						% Automatisxche Erstellung des Inhaltsverzeichnis, der Überschriften und Kopfzeilen
\usepackage{setspace}						% Zeilenabstand kann einsgestellt werden
\usepackage{ulem} 							%Durchstreichen von Text
\usepackage{longtable}						% Verwendung von mehrseitige Tabellen ermöglichen
\usepackage{multirow,bigdelim}
\usepackage{booktabs} 						% Besseraussehende Tabellen -> horizontale Linien mit \toprule, \midrule, \bottomrule
\usepackage{parcolumns}						% Zusammenfassen mehrerer Spalten in eiener Zeile
\usepackage{tabularx} 						% Gezielte Einstellungen der Breiten und Höhen von Tabellen
\usepackage{threeparttable}			  % Environment for creating a table with footnotes
\usepackage{graphicx}							% Einbindung von Abbildungen mit vielen Einstellmöglichkeiten
\usepackage{flafter}							% Abbildung werden niemals vor Verweis setzen
\usepackage{placeins}							% Setzt alle gesetzten Abbildungen vor dem Befehl \FloatBarrier
\usepackage{rotating}							% Erlaubt Abbildungen zu rotieren
\usepackage{tikz}								% Diagramme direkt in Latex erstellen
\usetikzlibrary{external}										% Plots als sperate PDF erstellen -> 						% Erlaubt die Einbindung von Diagrammen als Abbildung mit eigenem Tex-File
\usepackage{pgfplots}							% Diagramme direkt in Latex erstellen; Ergänmzung zu TIKZ
\usepackage{epstopdf}							% Hilfseinstellungen für Abbildungen bei PDF-Darstellung
\pgfplotsset{compat=newest}						% Einstellungen für TEX-Diagramme
\pgfplotsset{plot coordinates/math parser=false}% Einstellungen für TEX-Diagramme
\usepackage{pdfpages}							% Einbinden ganzer PDF-Seiten
\graphicspath{{../02_Grafiken/}} %Speicherort der Vektorgrafik
\usepackage{float}
\usepackage{mhchem}

\renewcommand{\min}[1][]{
	\ifthenelse{\isempty{#1}}{\operatorname{min}}{\ensuremath{\underset{#1}{\text{min}\,}}}
}

\def\LTT{
  Institute of Technical Thermodynamics, 
  RWTH Aachen University,
  Schinkelstraße 8,
  52062 Aachen,
  Germany
}

\newcommand{\mytitle}{Decision-Based vs. Distribution-Driven Clustering for \\ Stochastic Energy System Design Optimization}

\newcommand{\affil}{
  \begin{itemize}[leftmargin=3mm, itemsep=0mm]
    \item[$^a$]\LTT
  \end{itemize}
}

\def\firstAuthor{Boyung Jürgens}
\newcommand{\myauthor}{\firstAuthor$^{a}$, Hagen Seele$^{a}$, Hendrik Schricker$^{a}$, Christiane Reinert$^{a}$, and Niklas von der Assen$^{a,*}$}
\newcommand{\correspondingauthor}{N. von der Assen, \LTT \newline E-mail: niklas.vonderassen@ltt.rwth-aachen.de}
\author{\myauthor}

\usepackage[
  colorlinks,
  linkcolor=blue,
  citecolor=blue,
  urlcolor=blue,
  pdftitle={\mytitle},
  pdfauthor={\firstAuthor}
]{hyperref}
\usepackage[capitalise, nameinlink]{cleveref}
\crefname{table}{Tab.}{Tab.}

\newcommand{\setpgfexternalcounter}[1]{
  \makeatletter%
  \pgfkeysgetvalue{/tikz/external/figure name}\myexternalname
  \expandafter\gdef\csname c@tikzext@no@\myexternalname\endcsname{#1}%
  \makeatother
}

\begin{document}

  \thispagestyle{firststyle}

  \begin{center}
    \begin{large}
      \textbf{\mytitle}
    \end{large} \\
    \myauthor
  \end{center}

  \vspace{0.5cm}

  \begin{footnotesize}
    \affil
  \end{footnotesize}

  \vspace{0.5cm}

  \begin{abstract}
Stochastic programming is widely used for energy system design optimization under uncertainty but can exponentially increase the computational complexity with the number of scenarios. Common scenario reduction techniques, like moments-matching or distribution-driven clustering, pre-select representative scenarios based on input parameters. In contrast, decision-based clustering groups scenarios by similarity in resulting model decisions. Decision-based clustering has shown potential in network design and fleet planning. However, its potential in energy system design remains unexplored.

In our work, we examine the effectiveness of decision-based clustering in energy system design using a four-step method: 1) Determine the optimal design for each scenario; 2) Aggregate and normalize installed capacities as features reflecting optimal decisions; 3) Use these features for k-medoids clustering to identify representative scenarios; 4) Utilize these scenarios to optimize cost in stochastic programming.

We apply our method to a real-world industrial energy system modeled as a mixed-integer linear program. We incorporate uncertainty by scaling time series with representative factors. We generate 500 single-year scenarios via Monte Carlo sampling, which we reduce using decision-based clustering. For benchmarking, we conduct distribution-driven k-medoids clustering based on the representative factors. In our case studies, both clustering methods yield designs with similar cost efficiency, although decision-based clustering requires substantially more computational resources. To our knowledge, this is the first application of decision-based clustering on energy system design optimization. Future research should investigate the conditions under which decision-based clustering yields more cost-efficient designs compared to distribution-driven clustering.
\end{abstract}

\vspace{0.5cm}


\vspace{0.75cm}

\newpage
\pagestyle{laterstyle}
\section{Scenario Reduction in Optimization under Uncertainty}\label{sec:Intro}
Decarbonization across multiple sectors such as energy or industry requires a comprehensive transformation of energy systems. Energy system optimization modeling can support this transformation. However, as energy systems are characterized by long planning horizons, future uncertainty is a major challenge in energy system optimization modeling~\citep{Fodstad.2022}.

Stochastic programming (SP) can be used to incorporate uncertainty in optimization models. However, the computational effort required for SP can rise exponentially with the number of considered scenarios, potentially rendering the optimization intractable~\citep{Hewitt.2022}. Therefore, scenario reduction techniques are commonly employed to identify a smaller set of representative scenarios, for example through scenario clustering.

Scenarios are typically clustered using distribution-driven (DD) clustering, which groups them based on their input parameters~\citep{Mavromatidis.2018}. In contrast, decision-based (DB) clustering groups scenarios that lead to similar model decisions in deterministic optimizations, i.e., optimizations where only a single scenario is considered ~\citep{Hewitt.2022}. Thereby, DB clustering considers the optimization problem’s underlying structure~\citep{Zhuang.2024}. To illustrate the motivation behind DB clustering, consider two scenarios that, despite having significantly different input parameters, both result in similar model decisions. By grouping these scenarios, the representative scenario retains information about the optimal model decisions~\citep{Hewitt.2022}. Although DB clustering has shown promise in fields such as stochastic network design and fleet planning~\citep{Hewitt.2022}, it has not yet been applied to the design optimization of multi-energy systems.

In this work, we apply DB clustering to the design of multi-energy systems. By incorporating model decisions in the clustering, we aim to identify representative scenarios that yield more cost-efficient designs than designs derived by DD clustering.

\FloatBarrier
\section[Scenario Clustering for Stochastic Energy System Optimization]{Scenario Clustering for Stochastic Energy System Optimization}\label{sec:Method}
In our work, we use a two-stage stochastic program $g^\text{st,D}(\Omega)$ that minimizes the expected total annualized cost $TAC$, where each scenario $\omega\in\Omega$ has a probability $\pi_\omega$:
\begin{align}
g^\text{st,D}\left(\Omega\right)
= \min\limits_{x^\text{D},x^\text{O}}{TAC(x^\text{D},x^\text{O},\Omega)}
= \min_{x^\text{D},x^\text{O}}
\left(
CAPEX\left(x^\text{D}\right)
+\sum_{\omega\in\Omega}{\pi_\omega OPEX\left(x^\text{D},x_\omega^\text{O},\Theta_\omega\right)}
\right)
\end{align}
Here, design decisions $x^\text{D}$ must be made before uncertainty is revealed (here-and-now).
Operational decisions $x^\text{O}$ are made after uncertainty is revealed (wait-and-see).
We model uncertainty by scaling uncertain parameters with representative factors $\Theta_\omega=\left[\theta_{\omega,f}\right]_{f\in F}$, where $F$ is the set of parameters described by the representative factors.
The operational expenditure $OPEX$ is affected by uncertainty and depends on time series data of utility demand, purchase cost of energy carriers, and availability of renewable energy.
The capital expenditure $CAPEX$ comprises maintenance cost and annualized investment cost of installed capacities.
In this work, $CAPEX$ is not directly affected by uncertainty.
When considering only a single scenario $\omega$, $g^\text{st,D}(\Omega)$ simplifies to a deterministic design optimization
\begin{align}
g^\text{det,D}\left(\omega\right)=\min_{x^\text{D}, x^\text{O}}{TAC\left(x^\text{D}, x^\text{O},\omega\right)}
=\min_{x^\text{D}, x^\text{O}}
\left( CAPEX\left(x^\text{D}\right)
+OPEX\left(x^\text{D}, x^\text{O},\Theta_\omega\right) \right).\label{eq:g_det_D}
\end{align}

\subsection{Decision-based Clustering}\label{sec:decBasClustering}
We propose a decision-based clustering approach that groups scenarios whose optimal designs derived from deterministic optimizations have similar installed capacities. Given $N$ original scenarios $\omega\in\Omega$ with equal probability, where $\left|\Omega\right|=N$:
\begin{enumerate}
	\item Solve $g^\text{det,D}(\omega)$ for each scenario and obtain the installed capacities $\mathbf{c}_\omega=\left[c_{\omega,t}\right]_{t\in T}$ of all technologies $t\in T$. Each $c_{\omega,t}$ is an entry in $x^\text{D}$.
	\item Define a set of technology categories $I$ and subdivide $T$ into disjoint aggregated subsets $T_i^\text{agg}\subseteq T$ with $\bigcup_{i\in I} T_i^\text{agg}=T$. Sum up installed capacities within $T_i^\text{agg}$ to get aggregated installed capacities: $\mathbf{c}_\omega^\text{agg}=\left[c_{\omega,i}^\text{agg}\right]_{i\in I}$ with $c_{\omega,i}^\text{agg}=\ \sum_{t\in T_i^\text{agg}} c_{\omega,t}$. Normalize $\mathbf{c}_\omega^\text{agg}$ by scaling each feature (i.e., the aggregated installed capacity of each technology category) to a range between 0 and 1 using the maximum installed capacities of each technology. Use these normalized values as features for clustering.
	\item Cluster the scenarios into $K$ (where $K\ll N$) representative scenarios $\Omega^\text{K,DB}$ using k-medoids with Euclidian distances. Each cluster centroid, i.e., each representative scenario, has a probability equal to the sum of the probabilities of all scenarios within that cluster.
\end{enumerate}

\subsection[Distribution-driven Clustering]{Distribution-driven Clustering}\label{sec:distDrivCluctering}
For benchmarking we employ a distribution-driven clustering approach, grouping scenarios with similar inputs. Given $N$ original scenarios $\omega\in\Omega$ with equal probability, where $\left|\Omega\right|=N$:
\begin{enumerate}
	\item Normalize the representative factors $\Theta_\omega$ by scaling each feature (i.e., each representative factor~$\theta_{\omega,f}$) to a range between 0 and 1.
	\item Cluster $K$ representative scenarios $\Omega^\text{K,DD}$ using k-medoids, analogous to step 3 of DB clustering, but with the normalized representative factors as features.
\end{enumerate}
Unlike DB clustering, DD clustering does not require solving $g^\text{det,D}(\omega)$, resulting in substantially reduced computational effort.

\subsection{Decision-based vs. Distribution-driven Clustering}\label{sec:DBvsDDCluctering}
To compare DB and DD clustering, we use representative scenarios obtained by both approaches ($\Omega^\text{K,DB}$ and $\Omega^\text{K,DD}$) to solve $g^\text{st,D}(\Omega^{\text{K},i})$, yielding the designs $x^\text{D,DB}$ and $x^\text{D,DD}$. The cost efficiencies of both designs are evaluated by solving deterministic operational optimizations, expressed as
\begin{align}
	g^\text{det,O}\left(x^\text{D},\omega\right)=\min_{x^\text{O}}{TAC\left(x^\text{D},x^\text{O},	\omega\right)}
	=\min_{x^\text{O}}
	\left( CAPEX\left(x^\text{D}\right)
	+OPEX\left(x^\text{D},x^\text{O},\Theta_\omega\right) \right).\label{eq:g_det_O}
\end{align}
Here, $x^\text{D}$ can represent either $x^\text{D,DB}$ or $x^\text{D,DD}$ , evaluated for each $\omega\in\Omega$. Eq.~\ref{eq:g_det_O} emerges from Eq.~\ref{eq:g_det_D} when constraining $x^\text{D}$ to specified installed capacities.

For a given scenario $\omega$, the regret $\rho\left(x^\text{D},\omega\right)=g^\text{det,O}\left(x^\text{D},\omega\right)-g^\text{det,D}\left(\omega\right)$ measures the absolute cost difference of a given design $x^\text{D}$ to the cost-optimal design derived by deterministic design optimization. Finally, we evaluate the relative cost difference between the designs $x^\text{D,DB}$ and $x^\text{D,DD}$ for a given scenario~$\omega$ using the relative regret
\begin{align}
\tilde{\rho}\left(x^\text{D,DB},x^\text{D,DD},\omega\right)=\left(\rho\left(x^\text{D,DB},\omega\right)-\rho\left(x^\text{D,DD},\omega\right)\right)/\,\,\mathbb{E}_{\omega\in\Omega}\left(g^\text{det,D}\ \left(\omega\right)\right).
\end{align}

\newpage
\FloatBarrier
\section[Case Study]{Case Study: Sector-coupled Industrial Energy System}\label{sec:casestudy}
\begin{figure}[ht]
	\centering
	\includegraphics[width=1.0\textwidth]{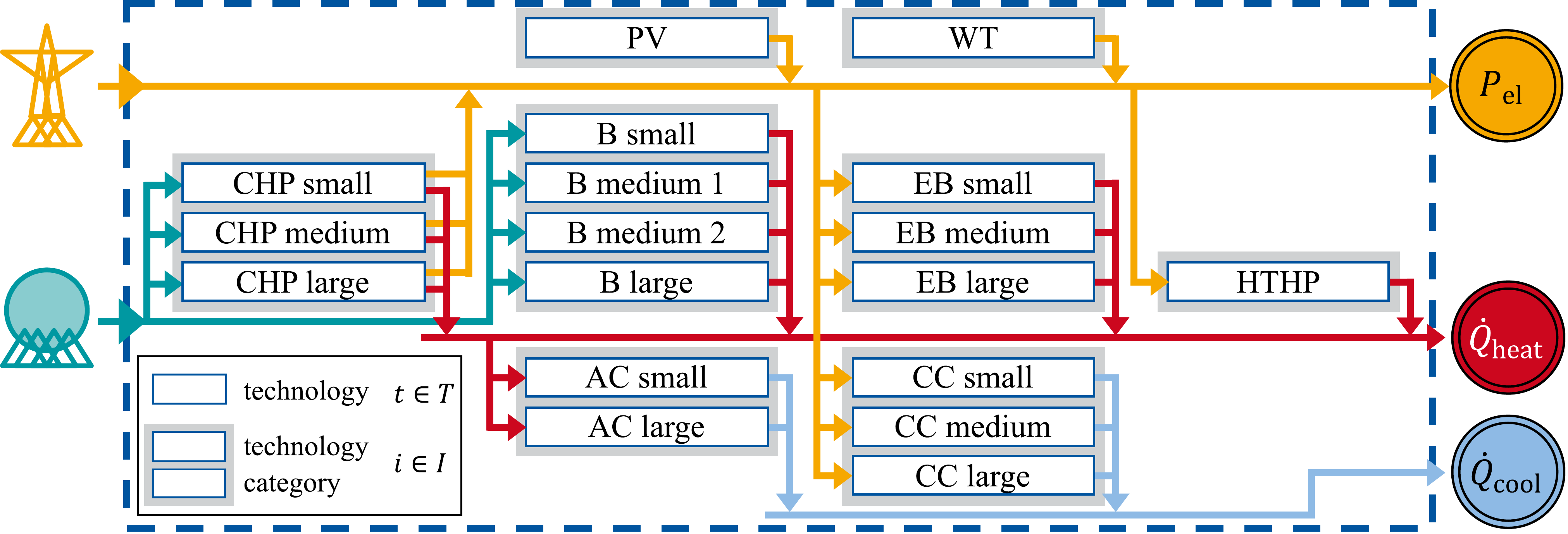}
	\caption{Superstructure of the industrial energy system serving electricity $P_\text{el}$, heating $\dot{Q}_\text{heat}$, and cooling $\dot{Q}_\text{cool}$ demand. The energy system is connected to the natural gas and electricity grid. Each component’s specifications can be found in~\cite{Reinert.2023}. The technology categories $I$ are: AC: absorption chiller, B: natural gas boiler, CC: compression chiller, CHP: combined heat and power, EB: electrode boiler, HTHP: high temperature heat pump, PV: photovoltaic, WT: wind turbine.}	
	\label{fig:case_study}
\end{figure}

We apply DB and DD clustering to reduce scenarios in a greenfield design optimization of an industrial energy system (Fig.~\ref{fig:case_study}). We model the energy system as a mixed-integer linear program accounting for minimal part-load and minimal installed capacities. We handle time series data with hourly resolution, aggregating the original 8760 timesteps into three typical days. We model uncertainty by scaling each time series using representative factors $\theta_{\omega,f}$, based on \cite{Seele.2023}: electricity $\theta_{\omega,\text{pe}}$ and natural gas $\theta_{\omega,\text{pg}}$ price; photovoltaic $\theta_{\omega,\text{ap}}$ and wind $\theta_{\omega,\text{aw}}$ availability; electricity $\theta_{\omega,\text{de}}$, heating $\theta_{\omega,\text{dh}}$, and cooling $\theta_{\omega,\text{dc}}$ demand. Additionally, we use combined representative factors: $\theta_{\omega,\text{at}}$ scales both availabilities, and $\theta_{\omega,\text{dt}}$ scales all utility demands simultaneously. We use model equations, component data, and time series data from \cite{Reinert.2023} \footnote{In contrast to~\cite{Reinert.2023}, we assume constant full-load efficiencies for part-load operation. In addition to the scaling factors applied by $\theta_{\omega,f}$, we scale the time series data for electricity prices ($\times0.28$), for natural gas prices ($\times0.15$), and for cooling demands $\left(\times0.43\right)$.}.

To explore the impact of modeling uncertainty at varying levels of detail, we conduct two case studies CS-4 and CS-7. CS-4 represents a low level of detail with representative factors $\Theta_\omega^\text{CS-4}=[\theta_{\omega,\text{pe}},\theta_{\omega,\text{pg}},\theta_{\omega,\text{at}},\theta_{\omega,\text{dt}}]$.
In CS-7, the representative factors are $\Theta_\omega^\text{CS-7}=[\theta_{\omega,\text{pe}},\theta_{\omega,\text{pg}},\theta_{\omega,\text{ap}},\theta_{\omega,\text{aw}},\theta_{\omega,\text{de}},\allowbreak\theta_{\omega,\text{dh}},\theta_{\omega,\text{dc}}]$.
Since CS-7 considers each representative factors individually, it offers a more detailed representation of uncertainty.

We use Monte Carlo sampling to generate $N=500$ scenarios for $\Omega$. Here, all representative factors are normally distributed with a mean of 1 and a standard deviation of 1/3. We allow only values of $\theta$ within $[0.5, 1.5]$. Scenarios violating this cutoff are resampled. We solve all optimizations to an optimality gap of $0.1\%$. 

\begin{figure}[ht]
	\centering
	\includegraphics[width=1.0\textwidth]{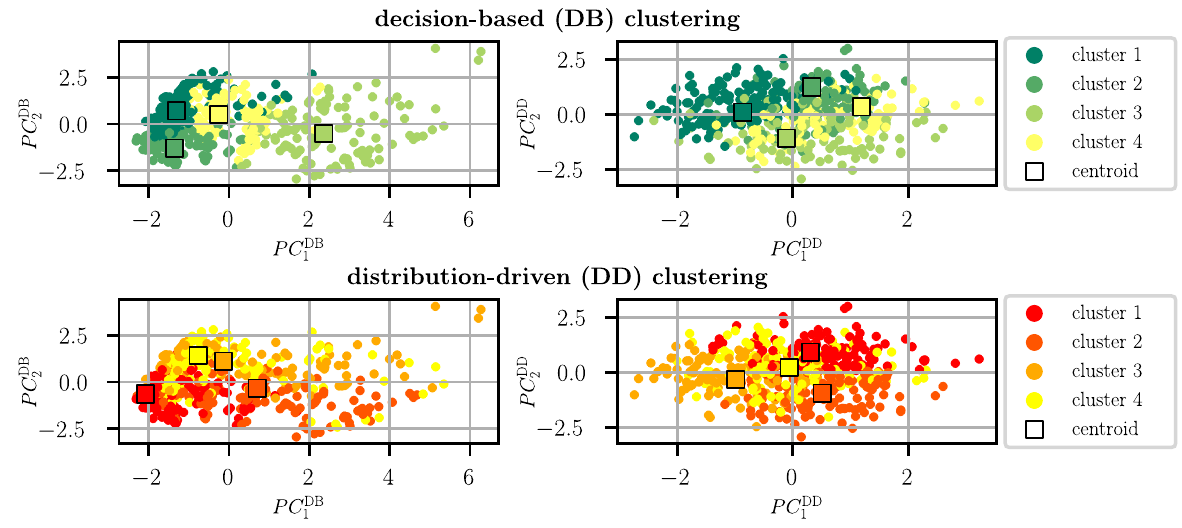}
	\caption{Principle component analysis (PCA) for CS-4. Top: decision-based (DB) clusters. Bottom: distribution-driven (DD) clusters. Left: PCA of the decision space (aggregated installed capacities ${\mathbf{c}_\omega^\text{agg}\ \forall\ \omega\in\Omega}$). Right: PCA of the input space (representative factors $\Theta_\omega^\text{CS-4}\ \forall\ \omega\in\Omega$).}	
	\label{fig:pca}
\end{figure}

We perform DB and DD clustering with the scenarios $\Omega$ from CS-4 and CS-7 to determine $K=4$ representative scenarios. DB clustering groups by eight features ($c_{\omega,i}^\text{agg}$ for each technology type in Fig.~\ref{fig:case_study}). DD clustering uses representative factors $\Theta_\omega$ as features. To visualize the clusters and centroids for both clustering approaches, we perform a principal component analysis on the cluster features. For CS-4, Fig.~\ref{fig:pca} maps DB and DD clusters onto the input space and decision space. Here, DB and DD clustering provide different clusters and centroids: While DB clusters form distinct areas in the decision space (Fig.~\ref{fig:pca}, upper left), DD clusters blend into each other in the decision space (Fig.~\ref{fig:pca}, bottom left). Thus, we conclude that clustering by representative factors does not effectively group scenarios with similar installed capacities in this case study.
\newpage
We use representative scenarios $\Omega^\text{K,DB}$ and $\Omega^\text{K,DD}$ to determine the optimal design by solving\linebreak $g^\text{st,D}(\Omega^{\text{K},i})$. To ensure feasibility in all original scenarios 
$\omega\in\Omega$, we include an extreme scenario with maximum utility demands in $g^\text{st,D}(\Omega^{\text{K},i})$ with a probability of zero. To determine the regrets $\rho(x^{\text{D},i},\omega)$ associated with the resulting designs $x^\text{D,DB}$ and $x^\text{D,DD}$, we compute $g^\text{det,O}(x^{\text{D},i},\omega)$ for each $\omega\in\Omega$. 

\begin{figure}[ht]
	\centering
	\includegraphics[width=1.0\textwidth]{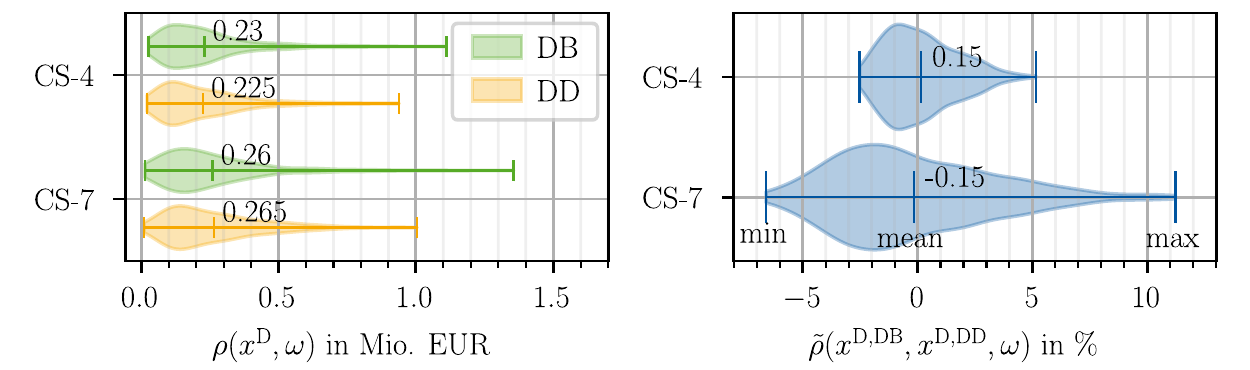}
	\caption{Distribution of regret $\rho\left(x^\text{D},\omega\right)$ (left) and relative regret $\tilde{\rho}(x^\text{D,DB},\ x^\text{D,DD},\omega)$ (right) for representative scenarios $\Omega^\text{K,DB}$ and $\Omega^\text{K,DD}$ determined by DB and DD clustering, respectively.}	
	\label{fig:regret_violin}
\end{figure}

Fig. \ref{fig:regret_violin} shows $\rho$ and $\tilde{\rho}$ for CS-4 and CS-7. Fig.~\ref{fig:regret_violin} shows that the mean relative regret $\tilde{\rho}$ is close to zero for both case studies ($\pm0.15\%$). Thus, the designs derived by DB and DD clustering do not significantly differ in cost efficiency.
\FloatBarrier
\newpage	
\section{Conclusions}\label{sec:conclusion}
In this work, we present a decision-based (DB) clustering approach for stochastic energy system design optimization that groups scenarios based on installed capacities. Initially, we carry out deterministic design optimizations for each scenario to determine the optimal installed capacities. We then aggregate installed capacities of similar technologies, using these aggregated values as features for k-medoids clustering.

We compare the performance of DB clustering to distribution-driven (DD) clustering in two case studies. In both case studies, DB and DD clustering provide designs with similar cost efficiency when tested against the original set of scenarios. This similar performance of both clustering methods can potentially be explained by our work’s uncertainty modeling: By using representative factors to scale times series data, DD clustering utilizes features carrying a high level of information. However, we assume that each time series maintains a consistent qualitative trajectory for each scenario, which may not reflect real-world variability and thus might be a debatable assumption.

This research serves as a proof-of-concept for applying DB clustering to energy system design optimization by grouping scenarios based on installed capacities. Future research should assess how the modeling of uncertainty impacts the performance of DB versus DD clustering and identify conditions where a superior cost efficiency of DB clustering’s designs justifies the increased computational demands of DB clustering.
\section*{Acknowledgements}
BJ gratefully acknowledges the financial support of the German Federal Ministry for Economic Affairs and Climate Action (Grand number: 03EI1084B).

\FloatBarrier

  \bibliographystyle{apalike}
  \renewcommand{\refname}{Bibliography}  
  \bibliography{literature.bib}

\end{document}